%% file: PCM.tex
\documentclass[3p]{elsarticle}

\usepackage{amsmath}
\usepackage{color}
\usepackage{hyperref}

\begin{document}

\begin{frontmatter}

\title{Point Cloud Movement For Fully Lagrangian Meshfree Methods}

\author[affil1,affil2]{Pratik Suchde\corref{mycorrespondingauthor}}
\ead{pratik.suchde@itwm.fraunhofer.de}
\author[affil1]{J\"org Kuhnert}

\cortext[mycorrespondingauthor]{Corresponding author}

\address[affil1]{Fraunhofer ITWM, 67663 Kaiserslautern, Germany}
\address[affil2]{Department of Mathematics, University of Kaiserslautern, 67663 Kaiserslautern, Germany}

\begin{abstract}
In Lagrangian meshfree methods, the underlying spatial discretization, referred to as a point cloud or a particle cloud, moves with the flow velocity. In this paper, we consider different numerical methods of performing this movement of points or particles. The movement is most commonly done by a first order method, which assumes the velocity to be constant within a time step. We show that this method is very inaccurate and that it introduces volume and mass conservation errors. We further propose new methods for the same which prescribe an additional ODE system that describes the characteristic velocity. Movement is then performed along this characteristic velocity. The first new way of moving points is an extension of mesh-based streamline tracing ideas to meshfree methods. In the second way, the movement is done based on the difference in approximated streamlines between two time levels, which approximates the pathlines in unsteady flow. Numerical comparisons show these method to be vastly superior to the conventionally used first order method. 
\\
\end{abstract}


\end{frontmatter}

\section{Introduction}

A moving Lagrangian framework is commonly used while modeling fluid flow. It often provides better approximations than the fixed Eulerian framework for flows with open free surfaces and multiphase flows with moving interfaces. The Lagrangian framework has the further advantage of avoiding the non-linear advection term, and often provides a more accurate depiction of transport phenomena. However, this comes at the cost of generally having a more restrictive time-step size control, and having the need to take special care for several aspects of conservation. 

In mesh-based methods, moving the mesh causes the additional disadvantage of mesh distortion. To avoid this distortion and the need to remesh, a large class of semi-Lagrangian and Arbitrary Lagrangian-Eulerian~(ALE) methods have been developed~(see, for example, \cite{ALE, semiLagrangian}). To improve conservation properties of mesh-based Lagrangian methods, the so-called ideas of trace-back and volume adjustments for mass conservation are often used  \cite{Arbogast2010,Iske2004}. The ideas include adjusting the volume of individual elements or cells based on their traced-back entities, and constructing upstream vertices and cells based on their corresponding downstream ones. 

Meshfree methods provide a more natural fit to Lagrangian frameworks than mesh-based methods. They use the numerical basis of a set of arbitrarily distributed nodes without any underlying mesh to connect them. These nodes could either be mass carrying particles or numerical points. Movement of this set of nodes, referred to as a point cloud, in a Lagrangian framework could also lead to distortion. However, 
point cloud distortion is easier to fix, as point clouds can easily be adapted locally, especially in meshfree methods that use numerical approximation points instead of mass carrying particles. 

To improve conservation properties, trace-back ideas used in mesh-based methods have also been generalized to Lagrangian meshfree methods \cite{Iske2007}. These methods involve adjusting the volume or mass of particles or the physical properties of approximation points appropriately \textit{after} their locations have been updated. In this paper, we consider the process of updating the point locations itself. We try to improve conservation properties by addressing the question of the optimum process of point movement, instead of adjusting physical quantities after movement is performed. 

For incompressible flow, inaccurate movement of points or particles results in errors in volume conservation. This results in the introduction of a numerical compressibility. This problem is especially relevant for applications with open free surfaces. To solve this and similar problems, various `artificial displacement' ideas have been introduced, especially in the context of the meshfree Smoothed Particle Hydrodynamics~(SPH). These involve performing an extra movement step in addition to the Lagrangian movement (or, equivalently, the addition of an extra term in usual Lagrangian movement step). This additional movement  is artificial in the sense that it is not based solely on the fluid velocity, as is the nature of the Lagrangian framework; rather it is based on improving different aspects of the numerical solution, such as to conserve total incompressibility of the simulated fluid \cite{Pahar2016} or to prevent local clustering of particles \cite{Khayyer2017,Shadloo2013}. Both physical and non-physical arguments have been used for the same. These methods are referred to under various names, such as particle shifting, artificial particle displacement, corrective displacement and particle regularization. In contrast to such methods, the methods presented in this paper improve the main Lagrangian movement step without the introduction of any additional artificial movement.

We begin the paper by showing that the most widely used first order approach for movement of points is extremely inaccurate, and propose two new methods for the same. The first new method is based on generalizing mesh-based streamline tracing. In this method, streamline velocities are approximated by an ODE system, and points are moved along these approximated streamlines. This method proves to give very good approximations for quasi-stationary flow problems. A further new method is developed which considers movement according to the change of these approximated streamlines between consecutive time levels. Numerical simulations show this method to give much better results for rapidly changing flow profiles.

To illustrate the use of these methods of point cloud movement, we use a meshfree Generalized Finite Difference Method~(GFDM) based on weighted least squares procedures like that in \cite{Jefferies2015,Kuhnert2014}. However, the same could also be applied to other meshfree methods. Since the methods discussed here can be applied to both numerical points in approximation point based meshfree methods and to particles in mass carrying particle based meshfree methods, the words `point' and `particle' are used interchangeably. 

The paper is organized as follows. In Section~\ref{sec:ExistingMethods} we introduce the notation used, and present the first order method which is commonly used to move point clouds in Lagrangian meshfree methods. In Section~\ref{sec:MeshedStreamlineTracing}, we talk about mesh-based streamline tracing, and discuss why direct extension of these ideas to meshfree methods is challenging. We then propose new methods for point cloud movement in Section~\ref{sec:NewMethods}. Numerical results on the application of different movement methods are shown for the Lagrangian advection equation in Section~\ref{sec:Advection}, and for the Lagrangian Navier--Stokes equations in Section~\ref{sec:INSE}. The paper is then concluded with a short discussion on the work in Section~\ref{sec:Conclusion}.

\section{Point Cloud Movement - First Order Method}
\label{sec:ExistingMethods}
We consider the time-integration of point locations while proceeding from time level $t^n$ to $t^{n+1}$, with $\Delta t = t^{n+1} - t^n$ being the time step size. Times in between the two time levels are referred to by $t^n + \tau$. Point cloud movement is done by integrating the equation
\begin{equation}
	\label{Eq:MovementODE}
	\frac{D\vec{x}}{Dt} = \vec{v}\,.
\end{equation}
Bracketed superscirpts are used to refer to the time level. Thus, point locations at $t^n$ and $t^{n+1}$ are referred to by $\vec{x}^{(n)}$ and $\vec{x}^{(n+1)}$ respectively; and velocities by $\vec{v}^{(n)}$ and $\vec{v}^{(n+1)}$ respectively. The closed-form movement is referred to as $\Delta \vec{x} = \vec{x}^{(n+1)} - \vec{x}^{(n)}$. Throughout this discussion, we only consider methods which decouple the movement step from the remaining PDEs being solved. Further, within a time step, we consider the movement to be done before solving the remaining PDEs. Thus, $\vec{v}^{(n+1)}$ is unknown during movement. If the movement process would be done after the solution of remaining PDEs is computed, $\vec{v}^{(n-1)}$ and $\vec{v}^{(n)}$ would be replaced with $\vec{v}^{(n)}$ and $\vec{v}^{(n+1)}$ respectively, throughout the paper.
\subsection{First Order Movement}
The most commonly used approach for moving point clouds is by assuming the velocity to be constant throughout the time step.
\begin{equation}
	\label{Eq:FirstOrderMove}
 \Delta \vec{x} = \vec{v}^{(n)}  \Delta t \,,
\end{equation}
which is performed at each numerical point or particle. Several variations of this form of movement are used in different meshfree methods. The most common is for each point to move with its own velocity $\Delta \vec{x}_i = \vec{v}^{(n)}_i  \Delta t\,,\; \forall i$. Some variants of SPH use movement with the average velocity in the neighbourhood of each particle \cite{Monaghan1989}, while others perform two first order movements per time-step: one based on an `intermediate' velocity, and one based on a `final' velocity \cite{Pahar2017}. 
Some meshfree methods use the average velocity of the current and the previous time step \cite{SkillenSPH}, and even refer to the same as a second order method. We club all of these methods under first order methods, since they assume the velocity to be constant throughout the time-step.

\begin{figure}[!htbp]
	\centering
	\includegraphics[width=0.2\textwidth]{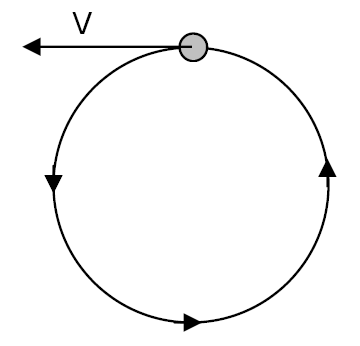}
	\caption{In flows with a rotational component, the most commonly used first order movement causes points to move along the tangential velocity.}
	\label{Fig:TangentialVelocity}%
\end{figure}

This method provides a very inaccurate approximation of moving points. The most significant errors come in capturing rotational parts of the flow correctly. This can be illustrated by the simple case of a rotating disc. Each point on the boundary of the disc has an instantaneous velocity in the tangential direction, as illustrated in Figure~\ref{Fig:TangentialVelocity}. Moving boundary points along this direction will result in the disc to constantly increase in volume. 

Despite the significant inaccuracies due to this method of movement, it is still used extensively today in almost all Lagrangian meshfree methods (see, for example, recent work in \cite{Luo2016,Pahar2017,SkillenSPH,VioleauSPH}).   

\section{Mesh-Based Particle and Streamline Tracing}
\label{sec:MeshedStreamlineTracing}

Particle tracing is a well established method for visualizations for fluid flow simulations \cite{Kipfer2003}. A common way to do this is to compute tangent curves by assuming a linearly varying velocity field \cite{Diachin1997,Nielson1999}. This amounts to computing analytical solutions to streamlines based on the underlying mesh. Such streamline tracing algorithms have also been extended to irregular geometries and more complex meshes \cite{Klausen2012}. 

These ideas of streamline tracing in visualization have also been extended to the movement of particles in Lagrangian particle methods. One such example is the Particle Finite Element Method~(PFEM) which uses moving Lagrangian particles with an underlying mesh. Like the aforementioned methods, vector fields are assumed to be piecewise linear on the background mesh. Based on these, piecewise integration of particle motion is performed by identifying the locations where particles cross simplex boundaries, like that done in \cite{Kipfer2003}. As a result, a closed-form solution to the movement is obtained \cite{Idelsohn2012,Nadukandi2015}. In addition to improving accuracy of movement, this method was shown to reduce the restriction on the time step size arising from the explicit nature of movement. Since most meshfree methods lack a global underlying mesh, such streamline tracing methods do not easily generalize to meshfree methods.

An alternative to analytical streamline tracing on the mesh is to approximate the streamlines numerically. This is done by splitting the time step into multiple sub-steps only for the purpose of determining the movement \cite{Idelsohn2012}. In each sub-step, each particle is moved with a first order method. Then, at the resulting temporary location at the end of a sub-step, a new velocity needs to be approximated for that temporary location which is then used for the movement in the next sub-step. This velocity approximation is also done based on the underlying mesh. Such a method has several drawbacks in a meshfree setting. Firstly, in the absence of a background mesh, the interpolation of velocities at the intermediate locations can prove to be very costly. Further, errors due to the interpolation at temporary locations can be quite significant \cite{Buning1988}. Moreover, since each sub-step uses first order movement of points, any improvement to the first order method can be coupled with this sub-stepping method.

Numerical integration based on high-order numerical methods such as Runge--Kutta methods have also been used for the purpose of particle tracking in mesh-based fluid flow problems. These methods have also been used for particle movement in molecular dynamics communities \cite{DEM}. However, these methods have been shown to be less accurate than those based on streamline tracing \cite{Lopes1998}. Moreover, as we shall show later, such higher order methods prove to be challenging in the context of many meshfree methods.

Due to these difficulties, the first order method, Eq.\,\eqref{Eq:FirstOrderMove}, is used almost exclusively for the purpose of updating particle locations. In this paper, we present methods to improve the accuracy of movement of Lagrangian points in meshfree methods. We begin by presenting a simple second order method and explaining why further higher order methods are not feasible. Then we extend streamline tracing to meshless contexts, and lastly, we combine the second order and streamline tracing methods to obtain a new method that gives better results for non-steady flows.

\section{Improved Methods for point cloud movement}
\label{sec:NewMethods}

\subsection{Second Order}
To avoid the inaccurate movement of the first order method, second order methods have been used by the present authors \cite{Jefferies2015,Kuhnert2014}. Such second order methods have also been used in ALE frameworks \cite{Hu2017_ALE} and in a select few SPH codes \cite{DualSPHysics}. Instead of assuming the velocity to be constant between two time levels, the velocity derivative is assumed to be constant. The movement can be considered to be done over a characteristic velocity $\vec{v}^c$, which satisfies the system
\begin{align}
	\frac{D\vec{v}^c}{Dt}(\tau) &= \frac{\vec{v}^{(n)} - \vec{v}^{(n-1)}}{\Delta t} \,, \;\phantom{abc}\;0<\tau<\Delta t\\
	\vec{v}^c(0) &= \vec{v}^{(n)}\,.
\end{align}
Note that this system is solved for each point. The particle displacements are then found by integrating along the characteristic velocity
\begin{equation}
	\label{Eq:MoveAlongCharacteristicVelocity}
	\frac{D\vec{x}}{Dt} = \vec{v}^c\,,
\end{equation}
to obtain
\begin{equation}
	\label{Eq:Move_M2}
	\Delta \vec{x} = \vec{v}^{(n)}\Delta t + \frac{1}{2}\frac{\vec{v}^{(n)} - \vec{v}^{(n-1)}}{\Delta t} (\Delta t)^2\,.
\end{equation}

We note that in this method, and henceforth, the variation of the characteristic velocity between two time steps is only considered for the sake of movement of points, and not for the computation of the new velocity, $\vec{v}^{(n+1)}$, which is done based on the relevant PDE being solved.

Higher order time integration methods, such as those used in Discrete Element Method~(DEM) simulations \cite{DEM}, are usually not possible, especially for approximation point-based meshfree methods. To prevent the distortion of point clouds, local adaptation is usually done. This involves adding points in locations containing `holes' and removing points in locations containing clusters of points \cite{Drumm2008}. All physical properties at the current time level are approximated at every new point created. However, approximating velocities from multiple time levels before the current one is extremely inaccurate. Thus, a newly created point would not posses $\vec{v}^{(n-2)}, \vec{v}^{(n-3)}, \dots$, which are needed in most higher order time integration methods. Thus, these higher order methods can not be used in many meshfree settings.

\subsection{Meshfree Movement Along the Streamline}

In fluid flow, streamlines describe the flow direction at a fixed instance in time. Thus, while going from time $t^{n}$ to $t^{n+1}$, moving points along the streamline would entail moving them along the velocity field at $t^{n}$. Meshed methods approximate this streamline velocity by performing linear interpolations on the underlying mesh. In contrast, we approximate the same by prescribing an ODE system which can be solved analytically. 

To perform movement along the velocity streamlines, we assume that each point is being advected based on the velocity gradient at its original location. Streamline velocities are computed based solely on the convective acceleration $(\vec{v}\cdot\nabla)\,\vec{v}$. Since the streamlines are taken at time $t^n$, the convective term can be taken as $(\nabla\vec{v}^{(n)})\,\vec{v}$, written as a matrix-vector product. Thus, the characteristic velocity between the time levels, along which movement is performed, can be taken by the initial value problem
\begin{align}
	\label{Eq:StreamlineODE}
	\frac{D\vec{v}^c}{Dt}(\tau) &= (\nabla\vec{v}^{(n)})\,\vec{v}^c\,, \;\phantom{abc}\;0<\tau<\Delta t\\
	\vec{v}^c(0) &= \vec{v}^{(n)}\,.
\end{align}
The above assumption of setting $\nabla\vec{v} \approx \nabla\vec{v}^{(n)}$ ensures that the resultant ODE system can be solved analytically to obtain
\begin{align}
	\vec{v}^c(\tau) &= \exp(\nabla\vec{v}^{(n)}\,\tau)\, \vec{v}^{(n)} \\&
	= \left[ \sum_{k=0}^\infty \frac{1}{k!} \left(\nabla\vec{v}^{(n)} \right)^k \tau ^k \right] \vec{v}^{(n)}\,.
\end{align}
Integrating the movement equation, Eq.\,\eqref{Eq:MoveAlongCharacteristicVelocity},
leads to the following closed form
\begin{equation}
	\label{Eq:Move_M3}
	\Delta\vec{x} = \left[ \sum_{k=0}^\infty \frac{1}{(k+1)!} \left(\nabla\vec{v}^{(n)} \right)^k \left(\Delta t \right) ^{k+1} \right] \vec{v}^{(n)}\,,
\end{equation}
%
where the velocity gradient is approximated numerically for each point based on the velocities of its neighbouring points. While this method accurately captures steady flow, a significant disadvantage is that it assumes the flow to be quasi-stationary. Ignoring the change of the velocity with time results in this method not providing very good approximations when there is a rapid change in the velocity profile.

\subsection{Movement According to the Change in Streamlines}

To overcome the disadvantage of the quasi-stationary flow assumption in the movement along the streamline method, we now present a method that is a generalization of the streamline method and the second order method. Rather than moving only along the streamline at the present time level, the difference of the streamlines of the present and the previous time levels is considered. For non-steady flows, this provides an approximation of the pathlines of the flow. The approximated streamline velocity at the present time level $t^n$ and previous time level $t^{n-1}$ are referred to as $\vec{v}^{s}$ and  $\vec{v}^{s_0}$ respectively. As done earlier, the streamline velocities are computed as
\begin{align}
&\left\{ \begin{aligned}
\frac{D\vec{v}^s}{Dt}(\tau) &= (\nabla\vec{v}^{(n)})\vec{v}^s\,, \;\phantom{abc}\;0<\tau<\Delta t\\
	\vec{v}^s(0) &= \vec{v}^{(n)}\,.
\end{aligned} \right.\\
&\left\{ \begin{aligned}
\frac{D\vec{v}^{s_0}}{Dt}(\tau) &= (\nabla\vec{v}^{(n-1)})\vec{v}^{s_0}\,, \;\phantom{abc}\;0<\tau<\Delta t\\
	\vec{v}^{s_0}(0) &= \vec{v}^{(n-1)}\,.
\end{aligned} \right.
\end{align}
Which leads to the streamline velocities
\begin{align}
\vec{v}^s(\tau) &= \left[ \sum_{k=0}^\infty \frac{1}{k!} \left(\nabla\vec{v}^{(n)} \right)^k \tau ^k \right] \vec{v}^{(n)} \,,\\
\vec{v}^{s_0}(\tau) &= \left[ \sum_{k=0}^\infty \frac{1}{k!} \left(\nabla\vec{v}^{(n-1)} \right)^k \tau ^k \right] \vec{v}^{(n-1)}\,.
\end{align}
Note that $\vec{v}^s $ is defined at time $t^n + \tau$, while $\vec{v}^{s_0}$ is defined at time $t^{n-1}+\tau$. The second order method assumed a characteristic velocity which had a constant derivative throughout the time step, based on the difference of velocities between the time levels. In this method, we assume a characteristic velocity such that the velocity derivative is based on the difference of the approximated streamline velocities.
\begin{align}
	\frac{D\vec{v}^c}{Dt}(\tau) &=  \frac{\vec{v}^s - \vec{v}^{s_0}}{\Delta t} \,, \;\phantom{abc}\;0<\tau<\Delta t\\
	\vec{v}^c(0) &= \vec{v}^{(n)}\,.
\end{align}
Integrating this leads to the characteristic velocity
\begin{equation}
	\vec{v}^{c}(\tau) = \vec{v}^{(n)} + \frac{1}{\Delta t}\left[ \sum_{k=0}^\infty \frac{\tau ^{k+1}}{(k+1)!} \left( \left(\nabla\vec{v}^{(n)} \right)^k \vec{v}^{(n)} - \left(\nabla\vec{v}^{(n-1)} \right)^k \vec{v}^{(n-1)} \right) \right] \,.
\end{equation}
Integrating again to obtain the displacement, 
\begin{equation}
	\label{Eq:Move_M4}
	\Delta\vec{x} = \vec{v}^{(n)}\Delta t + \frac{1}{\Delta t}\left[ \sum_{k=0}^\infty \frac{(\Delta t) ^{k+2}}{(k+2)!} \left( \left(\nabla\vec{v}^{(n)} \right)^k \vec{v}^{(n)} - \left(\nabla\vec{v}^{(n-1)} \right)^k \vec{v}^{(n-1)} \right) \right] \,.
\end{equation}
\begin{figure}
	\centering
	\def\svgwidth{0.4\columnwidth}
	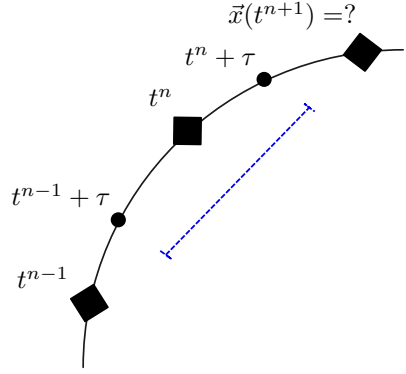
	\caption{Moving with the change of streamlines: The actual locations of the point being considered is represented by squares. The smaller circles represent instantaneous dummy locations between the time levels. The instantaneous derivative of the characteristic velocity is taken to be the difference in approximated velocities at the dummy locations~(the difference along the blue line). The movement is finally done along the resultant characteristic velocity.}
	\label{Fig:ChaneOfStreamline}%
\end{figure}
This process is represented in Figure~\ref{Fig:ChaneOfStreamline}. Unlike the streamline method, this method takes into consideration the change in the velocity field between time levels, and thus provides more accurate results for rapidly changing flows. This comes at the cost of the need to store $\nabla \vec{v}^{(n-1)}$.

In the numerical results presented later, the infinite summations in Eq.\,\eqref{Eq:Move_M3} and Eq.\,\eqref{Eq:Move_M4} are truncated after the first $5$ terms, which introduces a small error. Since each of the time-integration methods considered are explicit in nature, the difference in simulation times between the methods is not significant. The results are split into two parts. The first is a pure transport problem and the second incompressible flow according to the Navier--Stokes equations. 

\section{Advection Equation}
\label{sec:Advection}
We apply the different movement methods to a pure transport problem.

\begin{align}
	\frac{D\vec{x}}{Dt} &= \vec{v} \,,\\
	\frac{D\phi}{Dt} &= 0\,, \label{Eq:phiUpdate}
\end{align}
where $\vec{v}$ is the advection velocity, and $\phi$ is the physical quantity being transported. 

While both examples considered below use `simple' prescribed velocity fields, they are used to illustrate the impact of the point movement in Lagrangian transport in meshfree methods. In both cases, the prescribed velocity fields are such that the domain is undergoing rigid body motion. Thus, relative point positions do not change, and there is no deformation of the point cloud. Thus, there is no need to perform point additions or deletions to improve point cloud quality during the simulation, and errors in the needed interpolation do not affect the simulation. This, coupled with the fact that Eq.\,\eqref{Eq:phiUpdate} can be integrated trivially, ensures that the only source of error is the movement itself.

\subsection{Rotation}
\begin{figure}
  \centering
  \includegraphics[scale=0.4]{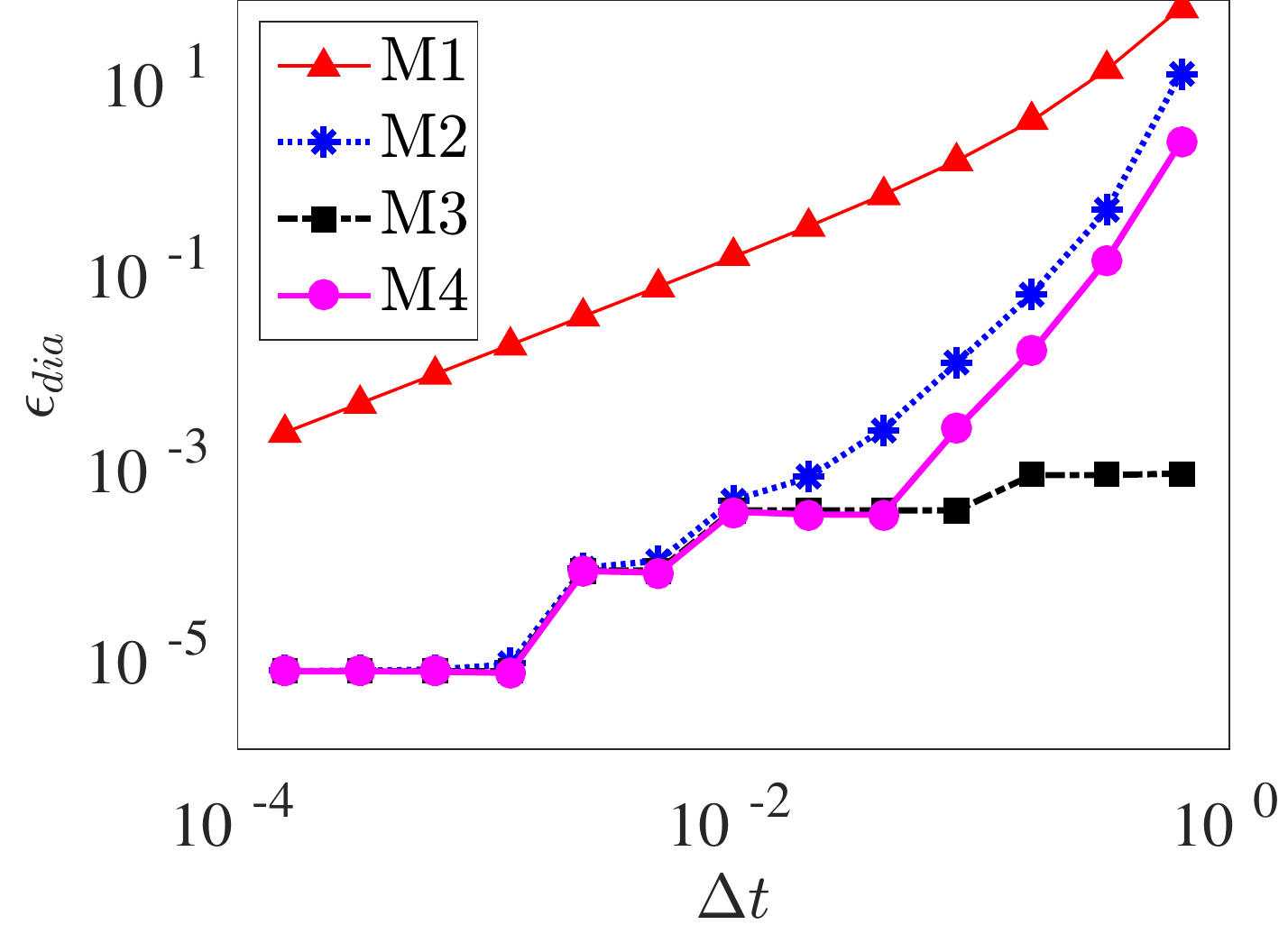}
  \caption{Errors for the rotating disc example. M1 stands for the first order method of movement, M2 the second order method, M3 the movement along the streamline, and M4 the movement according to the change of streamlines.}%
  \label{Fig:Rotation}
\end{figure}

We consider a circular disc of unit radius rotating about its center. The velocity field is given by
\begin{equation}
	\vec{v} = (-y,x)\,.
\end{equation}
The discretized disc contains $N=222$ points. Since only rigid body rotation is being performed, the area of the disc should be preserved. However, numerically, the disc expands due to movement of points along the tangential velocity at the boundaries, as shown in Figure~\ref{Fig:TangentialVelocity}. The error in the numerical result is measured as
\begin{equation}
	\label{Eq:EpsilonDia}
	\epsilon_{dia}= | d_{num} - d | \,,
\end{equation}
where $d_{num}$ is the numerical diameter of the disc at the end of two complete rotations, and $d$ is the theoretical diameter of the disc. The error with different time steps for the different movement methods is shown in Figure~\ref{Fig:Rotation}. The first order method produces the most inaccurate results. The remaining three methods produce similar results for small time steps, while the movement along the streamline produces the most accurate results for large time steps. For large time steps, the movement according to the change of streamlines is more accurate than the second order method. Even for the smallest time step considered, the first order method is worse than the other three methods by two orders of magnitude. 
\subsection{Transport along a Lissajous Curve}	

\begin{figure}
  \centering
  \includegraphics[width=0.45\textwidth]{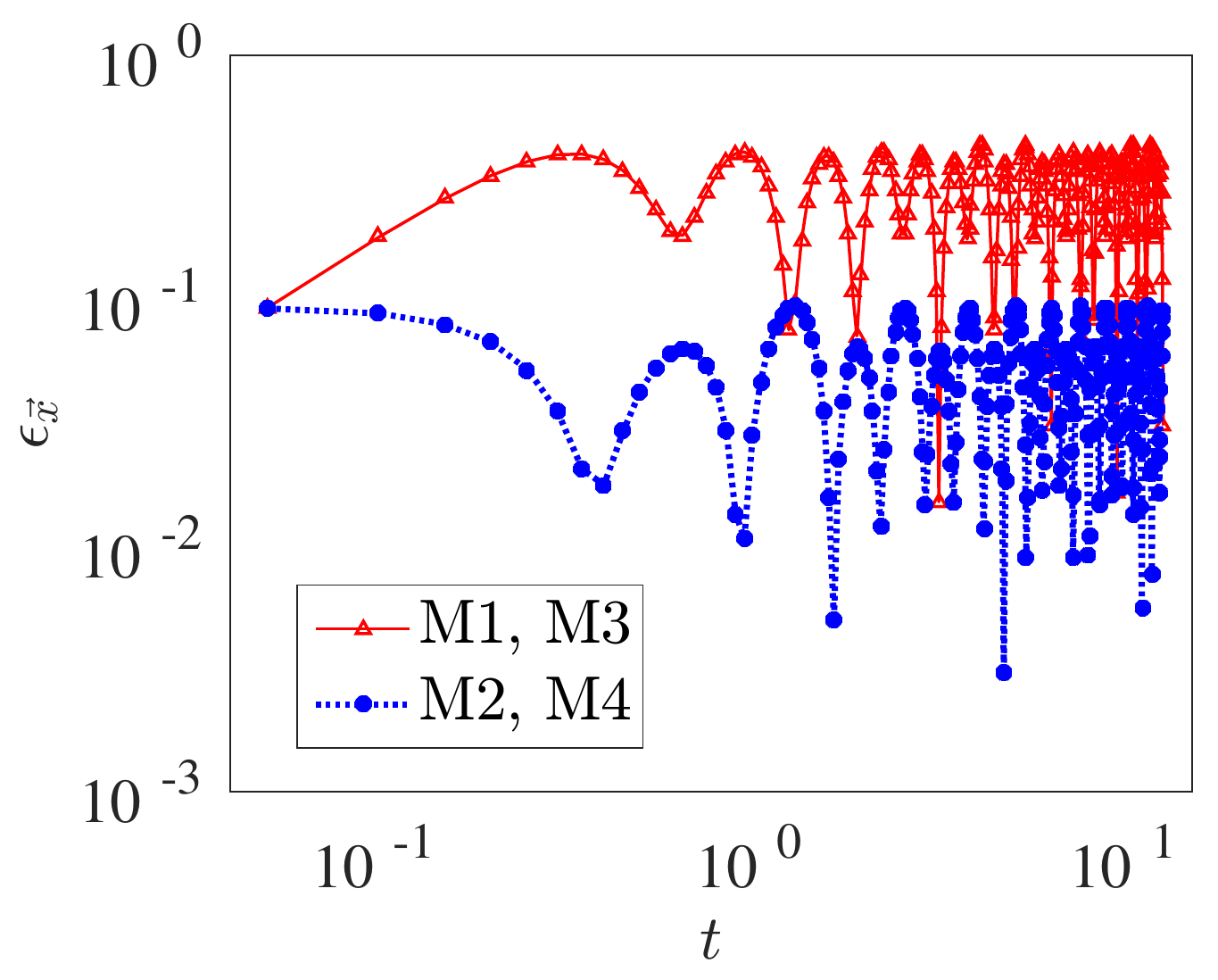}
  \caption{Errors in transport along Lissajous curves. M1 stands for the first order method of movement, M2 the second order method, M3 the movement along the streamline, and M4 the movement according to the change of streamlines.}%
  \label{Fig:Lissajous_results}
\end{figure}

We consider the same disc as in the previous example. The disc is being  transported along a time dependent velocity field
\begin{equation}
	\label{Eq:LissajousVelocity}   
	\vec{v} = ( 15\cos(5t+\frac{\pi}{2}) ,4\cos(4t))\,.
\end{equation}
For a disc with center starting at $(0,0)$, the exact trajectory of the center of the disc is given by the Lissajous curve
\begin{equation}
	\vec{x}_{exact} = ( 3\sin(5t+\frac{\pi}{2}) - 3 ,\sin(4t))\,.
\end{equation}
The spatially constant velocity means that $\nabla \vec{v} \equiv \mathbf{0}$, which results in Eq.\,\eqref{Eq:Move_M3} reducing to Eq.\,\eqref{Eq:FirstOrderMove}; and Eq.\,\eqref{Eq:Move_M4} reducing to Eq.\,\eqref{Eq:Move_M2}. Thus, in this case, the movement along the streamline and the first order method are the same; and the movement according to the change of the streamlines and the second order method are the same. Such cases with spatially constant velocity in a neighbourhood are relevant for flow regimes in the `far field' of the domain where the velocity does not vary much. The error in the numerical solution is measured by
\begin{equation}
	\epsilon_{\vec{x}} = \|\vec{x}_{num} - \vec{x}_{exact}\|\,,
\end{equation}
where $\vec{x}_{num}$ is the numerical center of the disc. For $\Delta t = 0.05$, the evolution of errors for the different methods of movement are shown in Figure~\ref{Fig:Lissajous_results}. The second order and change of streamline method produce more accurate results than the first order and streamline method.

\section{Incompressible Navier--Stokes Equations}
\label{sec:INSE}

Consider the incompressible Navier--Stokes equations in Lagrangian form
\begin{align}
	\frac{D\vec{x}}{Dt} &=  \vec{v} \,, \\
	\nabla\cdot\vec{v}  &= 0 \,,\\
	\frac{D\vec{v}}{D t} &= \frac{\eta}{\rho}\Delta \vec{v} - \frac{1}{\rho}\nabla p + \vec{g} \,,
\end{align}
where $\vec{v}$ is the fluid velocity, $p$ is the pressure, $\rho$ is the density, $\eta$ is the dynamic viscosity and $\vec{g}$ includes both gravitational acceleration and body forces. The numerical scheme consists of movement of the point cloud according to one of the methods mentioned earlier, which is followed by a meshfree projection method similar to the one in \cite{Jefferies2015,Kuhnert2014}. In addition, a $k-\epsilon$ turbulence model is also used. The boundary conditions used, including the ones at the free surface, are as detailed in \cite{Jefferies2015,Kuhnert2014}. While an inaccurate update of point locations represent one source of error in meshfree solution schemes to the incompressible Navier--Stokes equations, other sources of errors in the same have been discussed in our earlier work \cite{Suchde2017_CCC,Suchde2017_INSE}.

Inaccurate movement of the points on the free surface can result in a volume loss or gain, similar to that illustrated in the rotating disc example considered earlier. For incompressible flow with density $\rho$ fixed and constant throughout the domain, this error in volume conservation also represents an error in mass conservation. This error can be measured by the relative change in the total volume occupied by all points in the computational domain.
\begin{equation}
	\label{Eq:VolumeLoss}
	\epsilon_{V}=\frac{|\int_{\Omega_0} dV - \int_{\Omega_{end}} dV|}	{\int_{\Omega_0} dV},
\end{equation}
where $\Omega_0$ is the initial domain and $\Omega_{end}$ is the domain at $t_{end}$.

\subsection{Sloshing}
\begin{figure}
  \centering
  \includegraphics[width=0.35\textwidth]{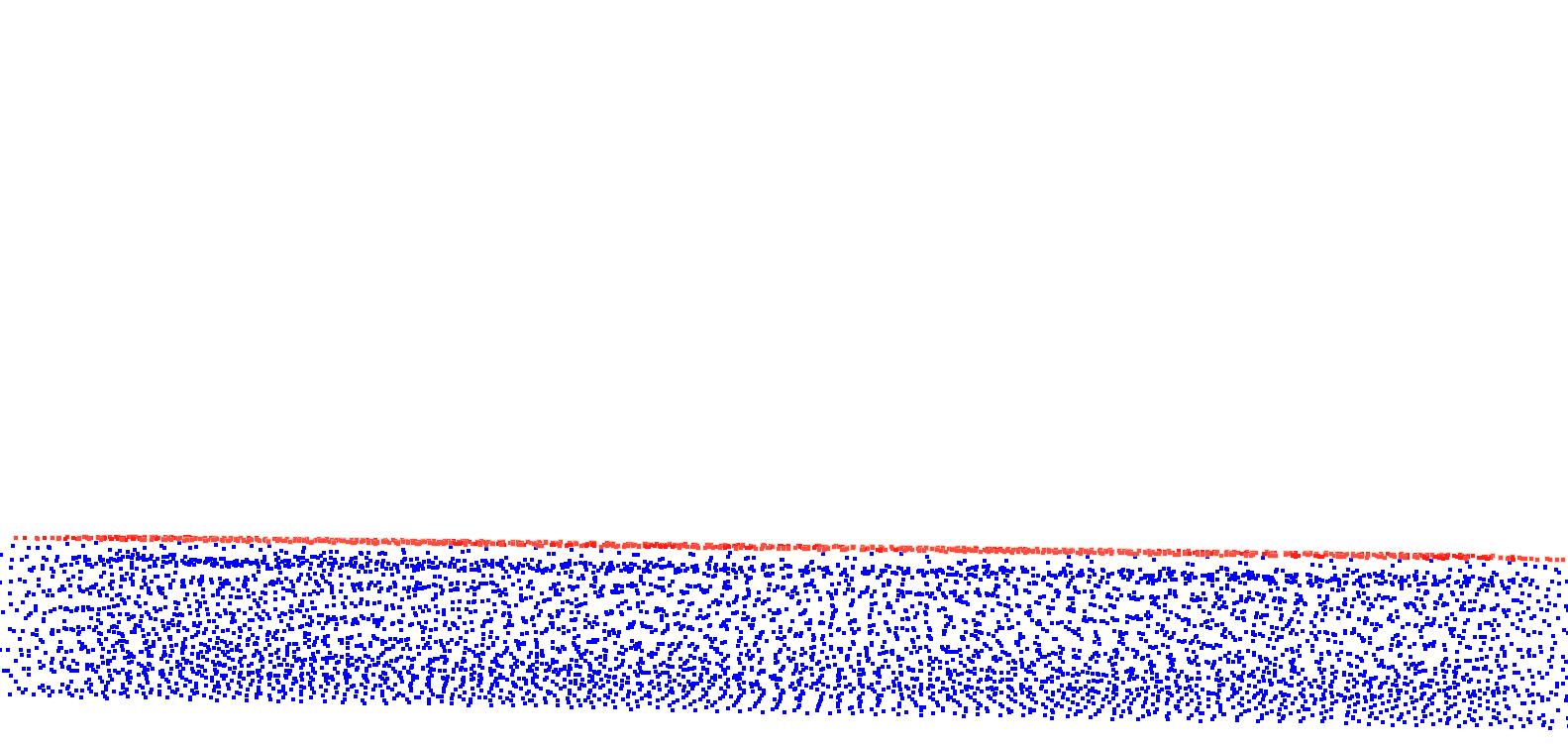}
  \includegraphics[width=0.35\textwidth]{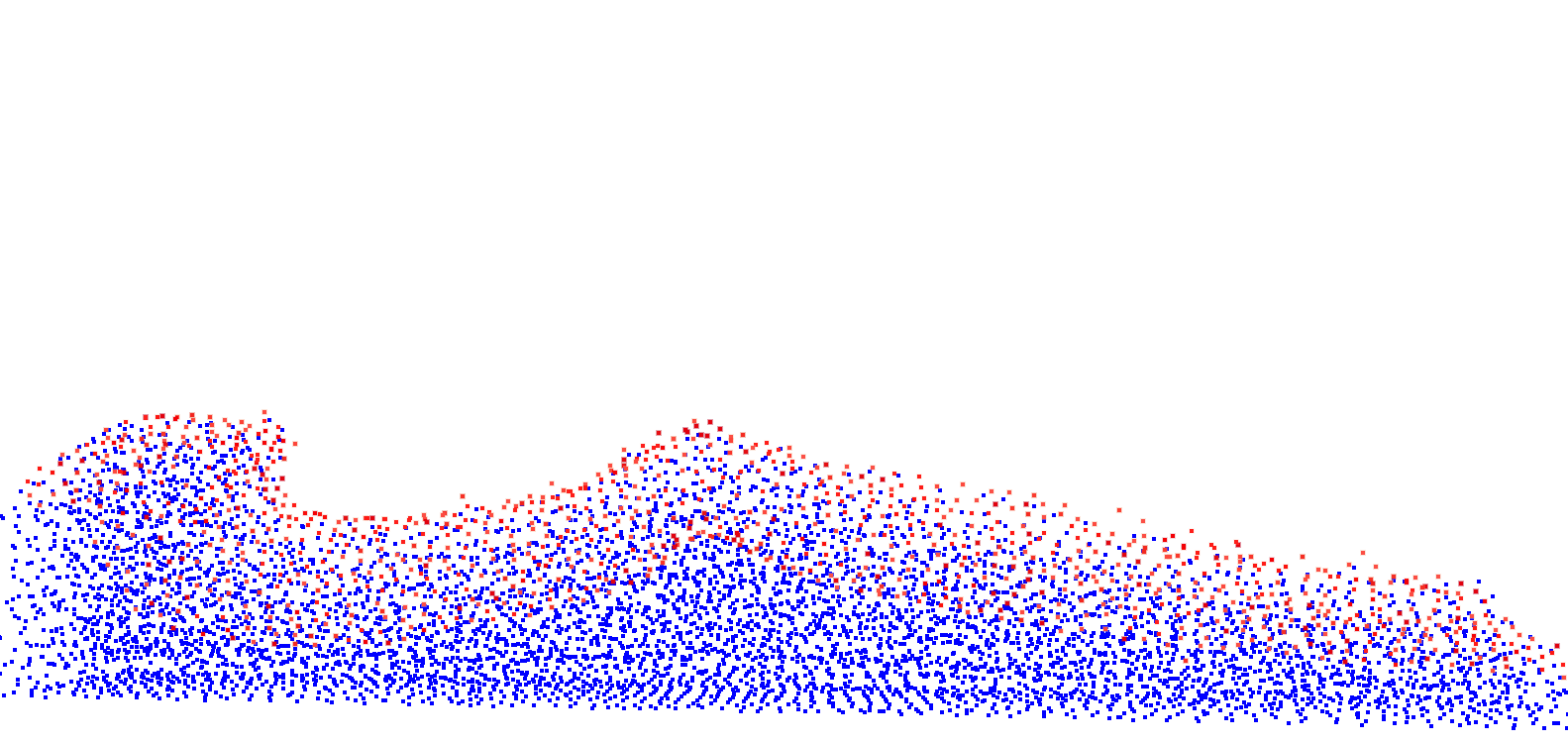}
  %
  \includegraphics[width=0.35\textwidth]{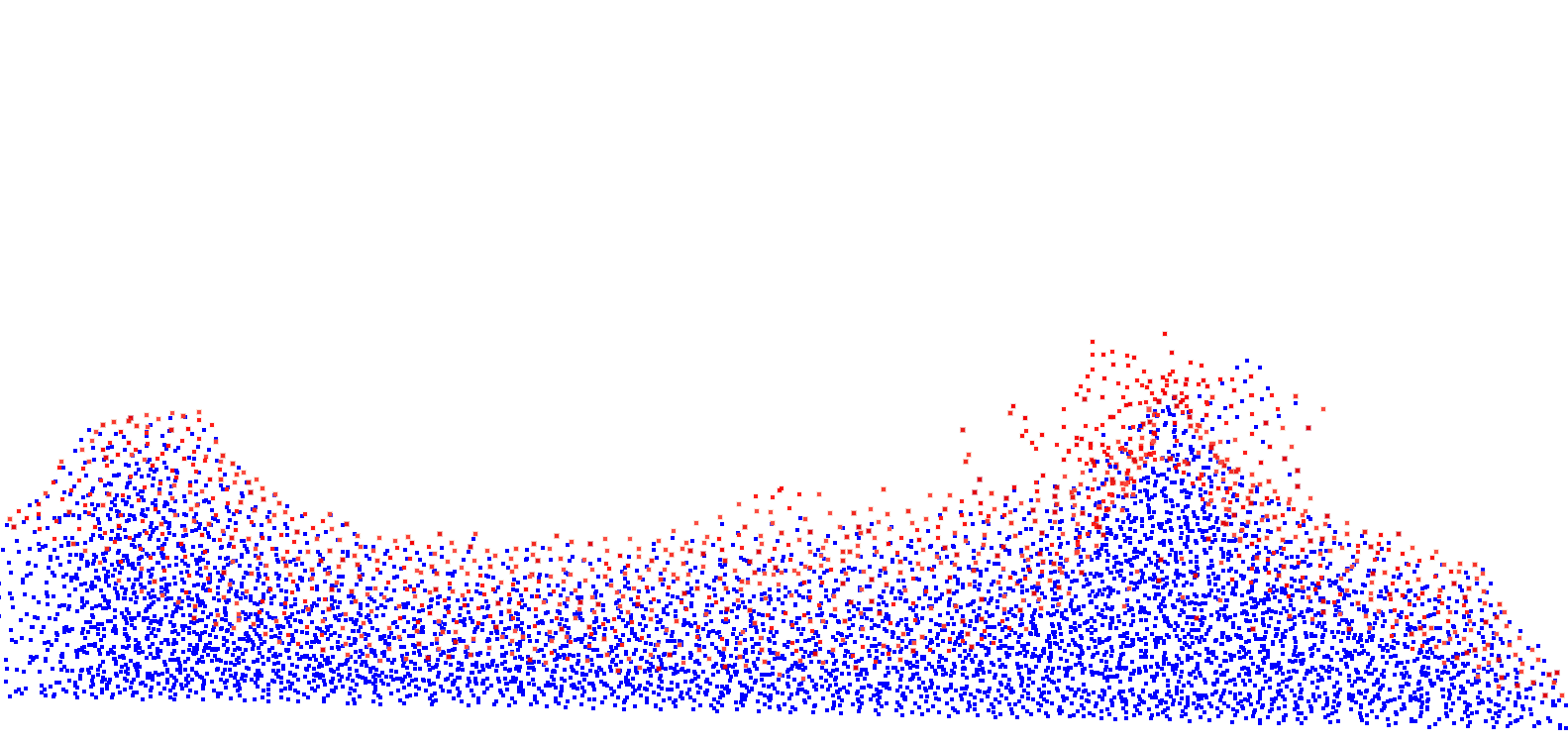}
  \includegraphics[width=0.35\textwidth]{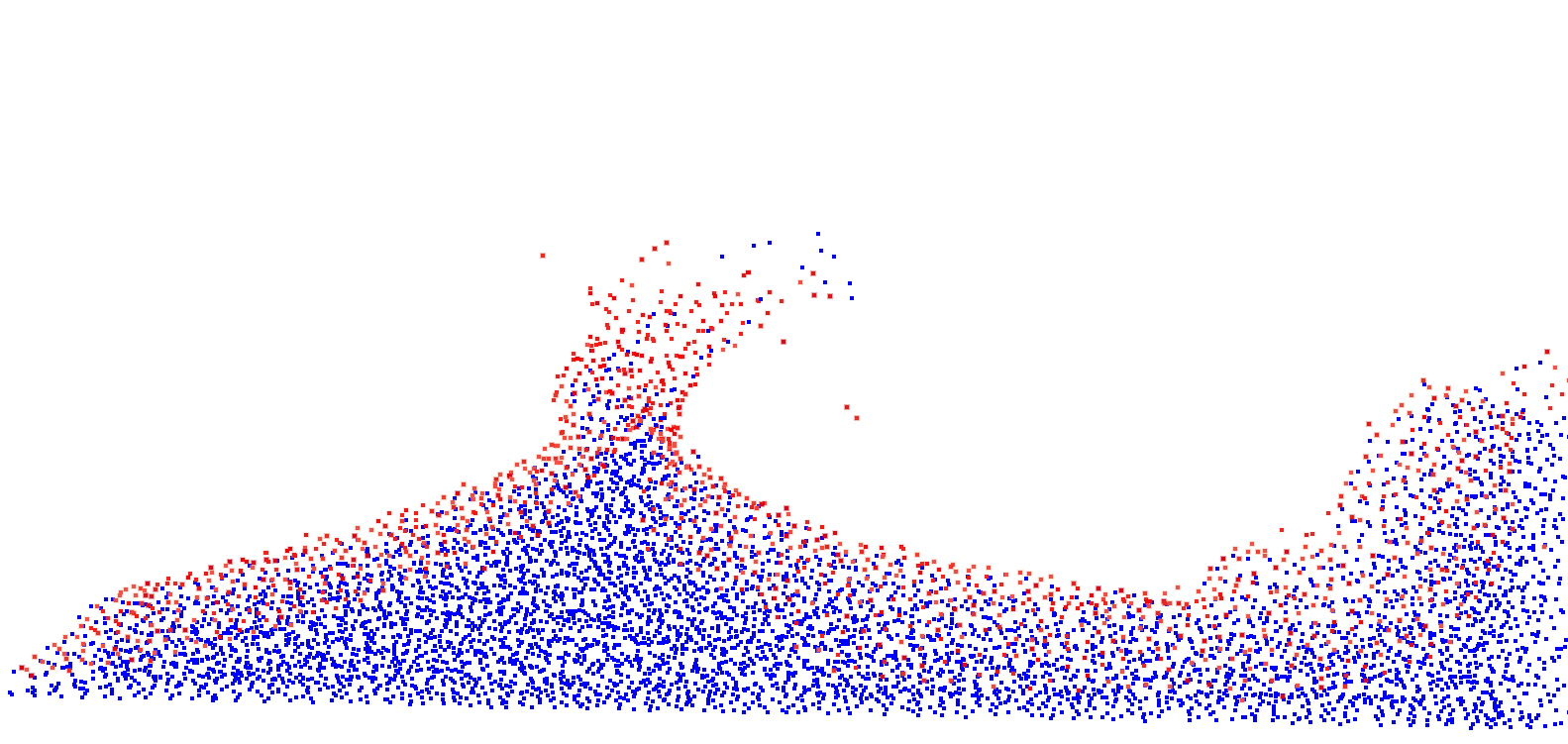}
  \caption{Sloshing. $t = 0s$~(top left), $t = 0.72s$~(top right), $t = 1.44s$~(bottom left) and $t = 1.92s$~(bottom right). The red color indicates the free surface.}
  \label{Fig:Sloshing3D}%
\end{figure}

We consider the test case of sloshing of water contained in a constantly moving rectangular box. This results in a non-stationary free surface as shown in Fig~\ref{Fig:Sloshing3D}. The initial state is taken to be at rest. Slip boundary conditions are used at the walls for the velocity. Homogeneous Neumann boundary conditions are used for the pressure at the walls. The movement of the box is represented in the gravitational and body forces term by setting $\vec{g}=(4\cos(7t),-10,-5)$. The simulation parameters are set as $t_{end} = 3s$, $\rho = 10^3 kg/m^3$ and  $\eta = 10^{-4} Pa\,s$, which results in a Reynolds number of the order of $10^7$. 
\begin{figure}[ht]
  \centering
  \includegraphics[scale=0.4]{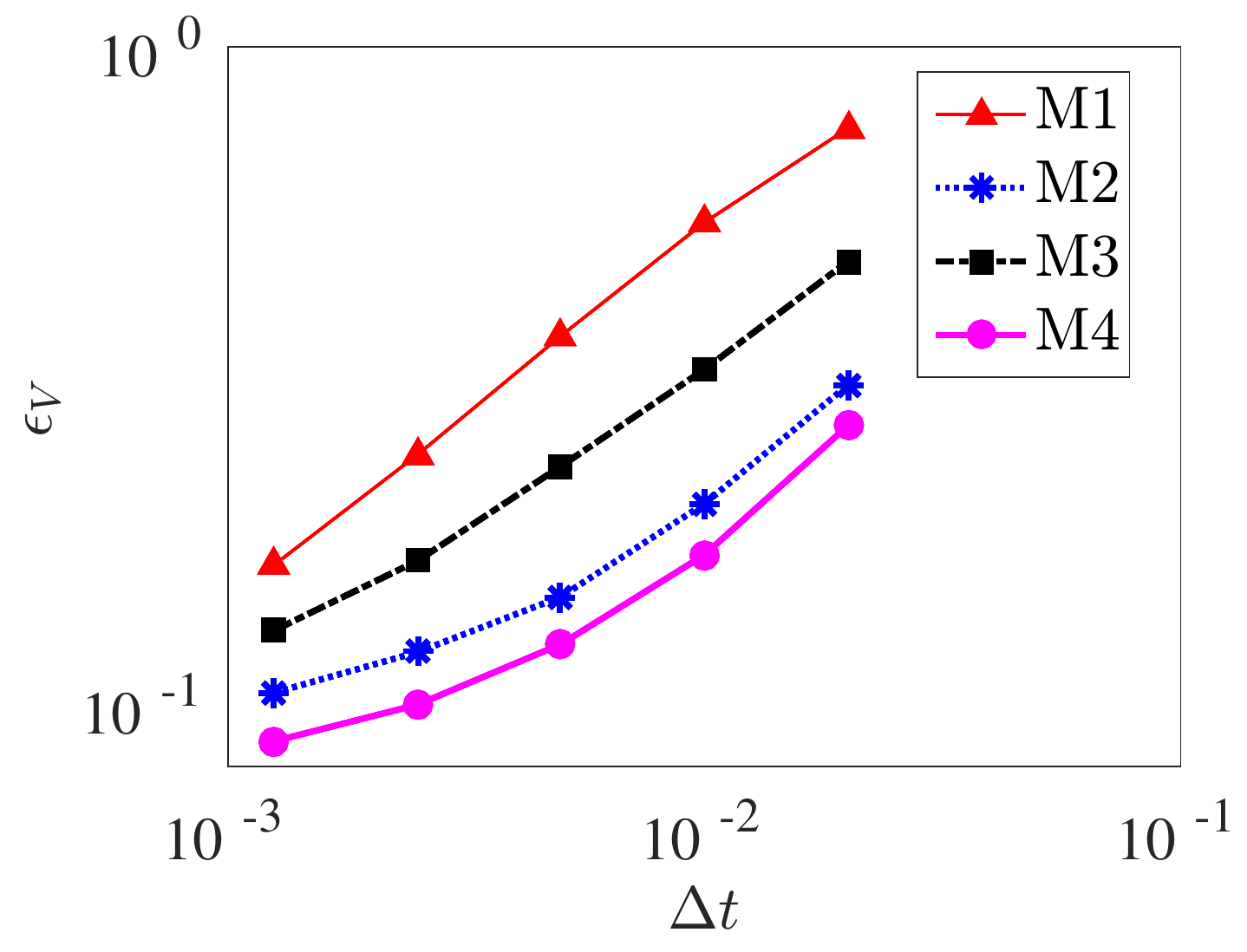}
  \caption{Volume conservation errors for the 3D sloshing test case. M1 stands for the first order method of movement, M2 the second order method, M3 the movement along the streamline, and M4 the movement according to the change of streamlines.}%
  \label{Fig:Sloshing3D_results}
\end{figure}

For a constant smoothing length of $h = 0.06$, which corresponds to an initial number of points $N = 12010$, Figure~\ref{Fig:Sloshing3D_results} shows the error in volume conservation for different time steps. For rapidly changing flows like that in this example, the movement along the streamline method is no longer the most accurate, which agrees with the theoretical expectation since this method does not take into account the velocity change with time. In comparison, the second order movement produces more accurate results. Further, point movement according to the change in streamline method produces slightly better results than the second order method.
\subsection{Tank Filling}

\begin{figure}
  \centering
  \includegraphics[width=0.3\textwidth]{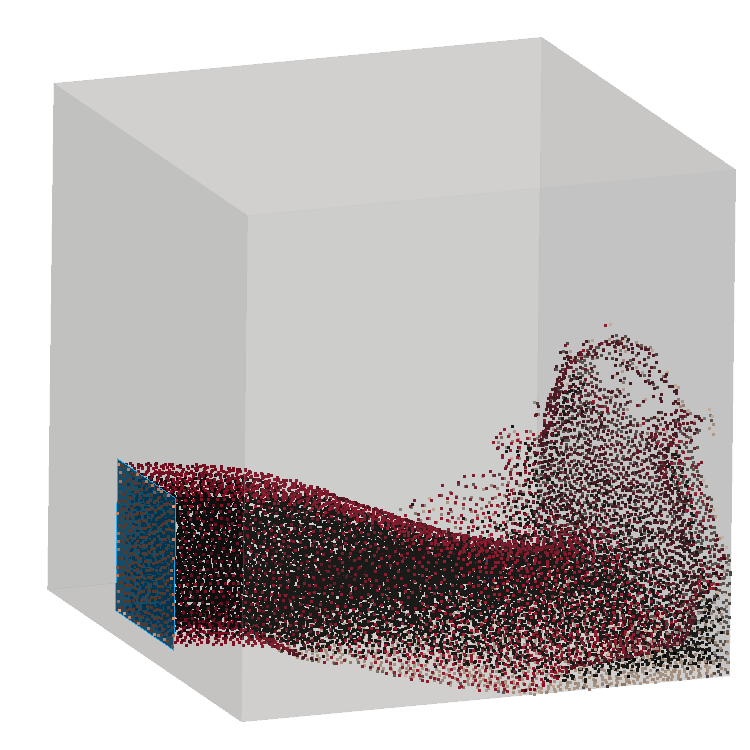}
  \includegraphics[width=0.3\textwidth]{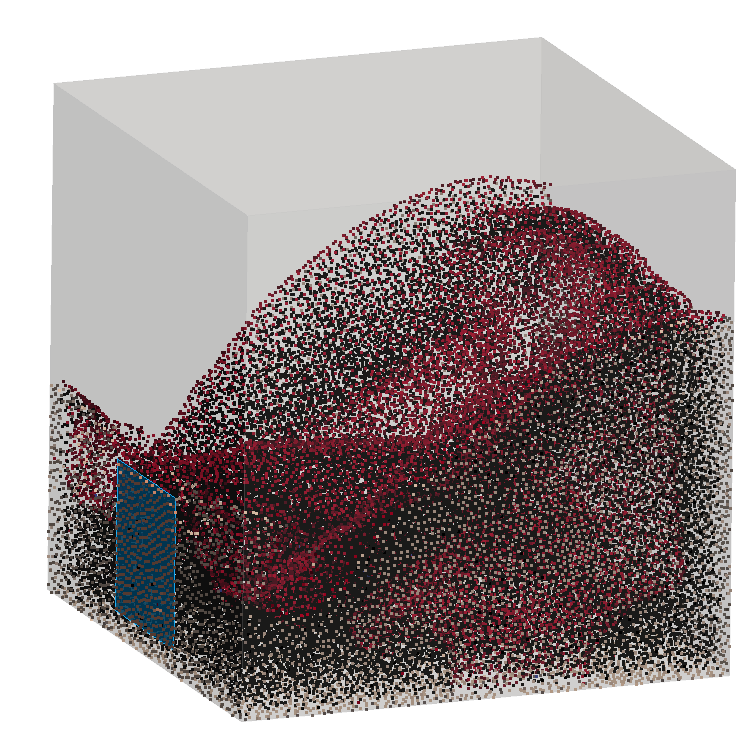}  
  \includegraphics[width=0.3\textwidth]{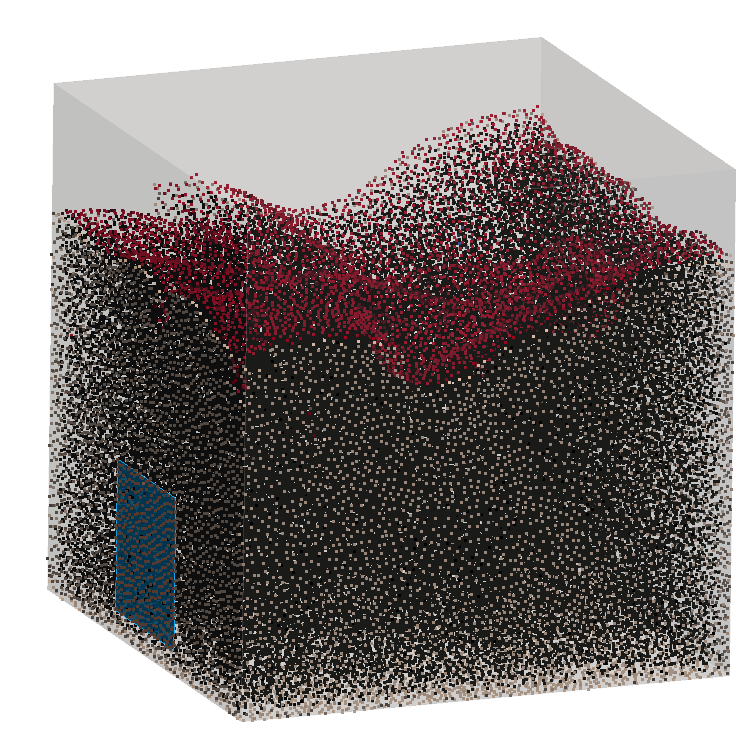}
  %
  %
  \caption{Tank filling. $t = 0.36s$~(left), $t = 1.32s$~(center), $t = 3.0s$~(right). Black points represent interior points, red points represent points on the free surface. The wall boundary points in white are hidden, except near the edges. The blue region on the left of the box marks the inflow. }
  \label{Fig:TankFillingImages}%
\end{figure}

We consider the filling of an initially empty cuboidal tank from a inlet near the bottom, as shown in Figure~\ref{Fig:TankFillingImages}. For a fixed inflow velocity, the theoretical total volume of the fluid inside the tank is known, and the numerical volume occupied by the point cloud can be compared with it.

This test case captures several challenges of fully Lagrangian meshfree methods. The impact of the high velocity jet on the domain walls can cause points to cluster near the location of impact. Incorrect movement of points near the impact further enhances this problem. Such clustering is a source of tensile instability in meshfree methods that use mass carrying particles, such as SPH. While mass conservation is seemingly trivial in these particle based meshfree methods, volume conservation is an issue that affects all meshfree Lagrangian methods. In all meshfree methods, the total geometric volume occupied by the point cloud or the particle cloud represents the physical volume occupied by the simulated fluid. Since incompressible flow is being considered, this volume should be conserved. However, this is not the case numerically, and thus, an artificial compressibility is introduced. Further, incorrect movement at and near the wall boundaries can lead to points escaping the simulation domain. Depending on the meshfree method being used, these escaped points are either projected back to the wall, or are deleted. Both these cases are another source of error in conservation.

\begin{figure}
  \centering
  \includegraphics[scale=0.4]{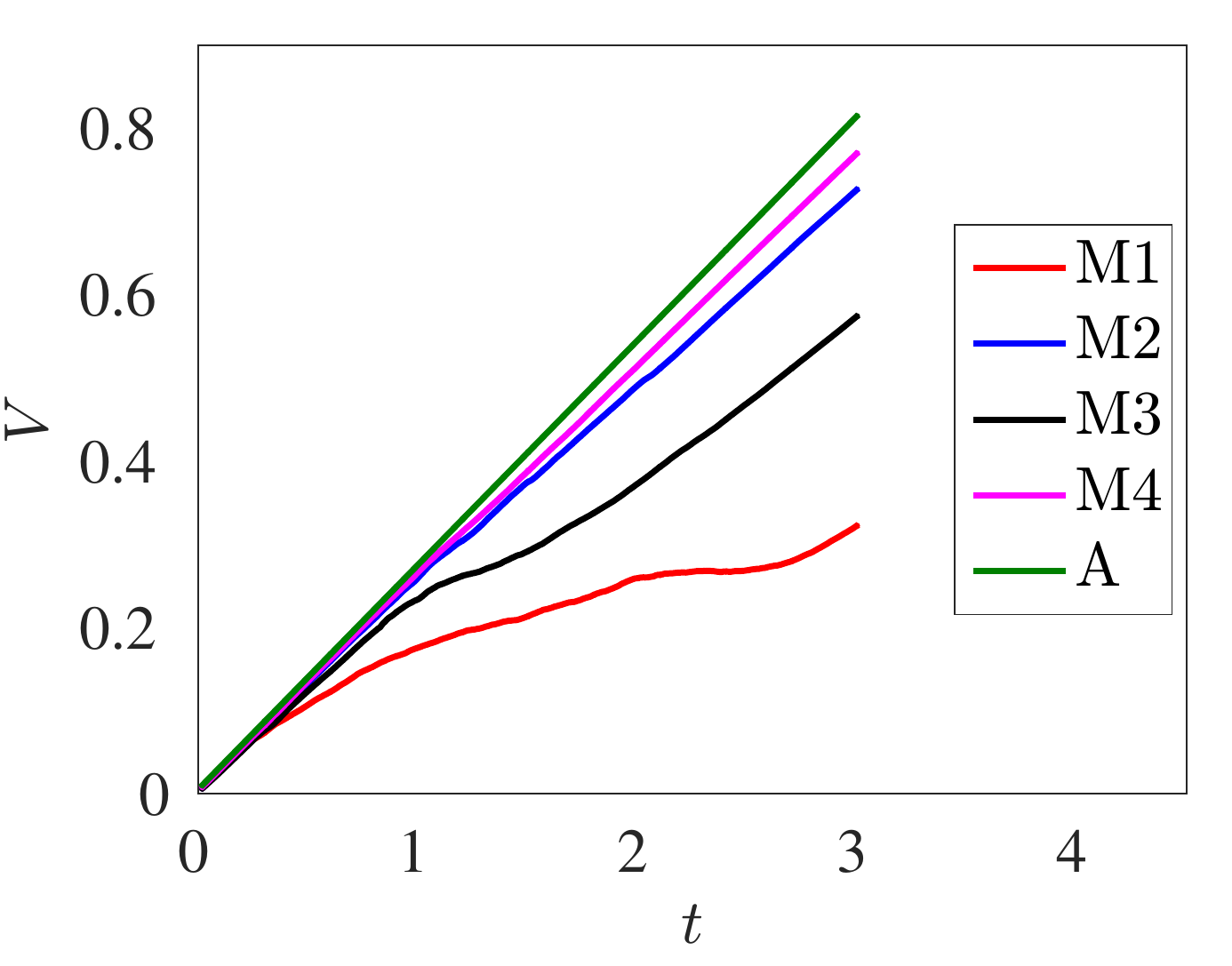}
  \includegraphics[scale=0.4]{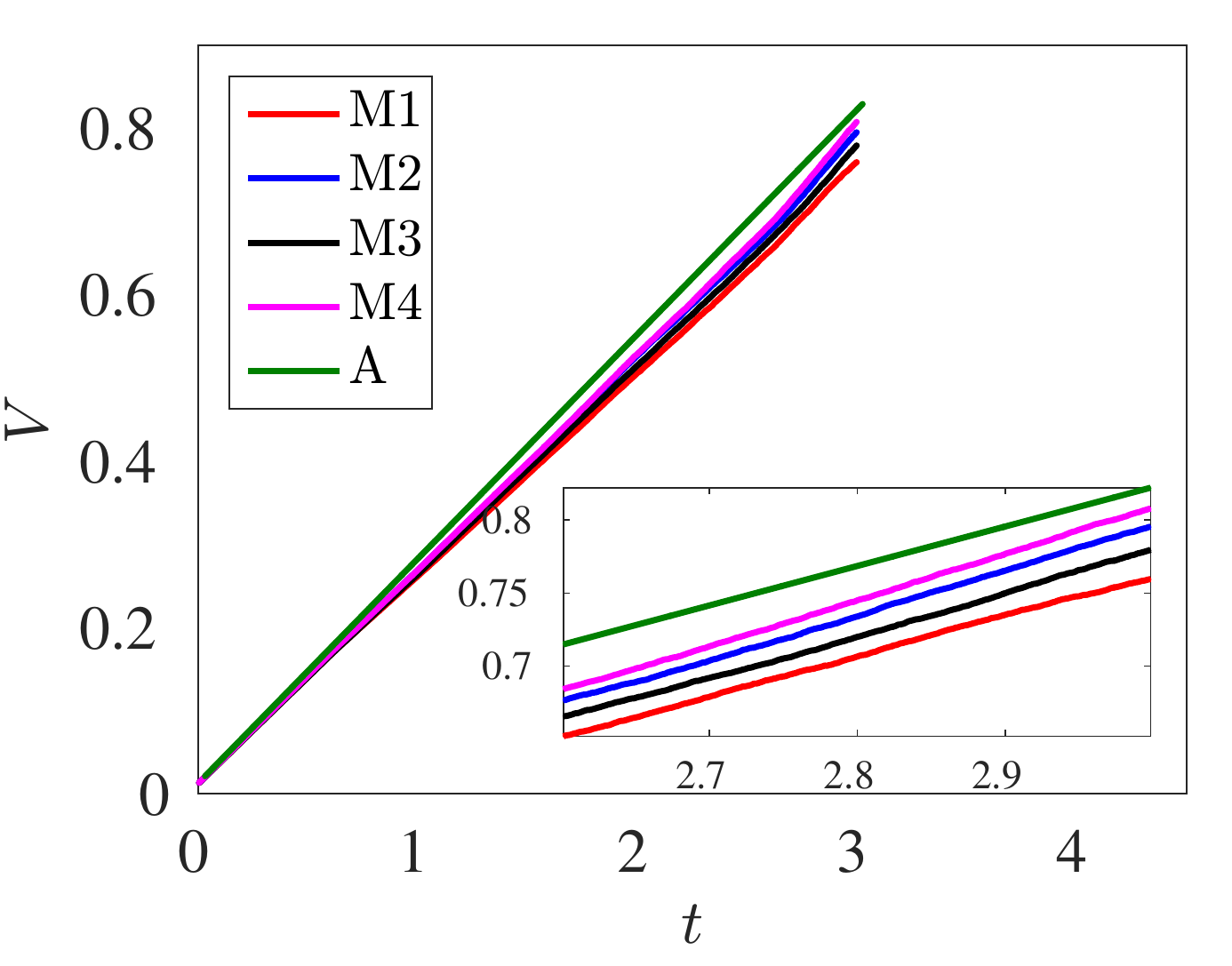}
  \caption{Total numerical volume for the tank filling test case. For a large time step of  $\Delta t = 0.6\times 10^{-2}$~(left), and a smaller one of  $\Delta t = 0.6\times 10^{-3}$~(right). M1 stands for the first order method of movement, M2 the second order method, M3 the movement along the streamline, and M4 the movement according to the change of streamlines. A represents the analytical volume.}%
  \label{Fig:TankFillingResults}
\end{figure}

The dimensions of the tank considered in the simulations are $1m \times 1m \times 1m$, and that of the inlet are $0.3m \times 0.3m$. A constant velocity inflow of $\vec{v}_{in}=(3,0,0)m/s$ is used at the inlet. Slip conditions are used on the walls. The inflowing fluid has properties of $\rho = 10^3 kg/m^3$ and $\eta = 10^{-3} Pa\,s$ which results in a Reynolds number of the order of $10^6$. The evolution of the numerical volume for simulations with all four movement methods considered are shown in Figure~\ref{Fig:TankFillingResults} for a large time step of $\Delta t = 0.6\times 10^{-2}$, and a small time step of $\Delta t = 0.6\times 10^{-3}$. In both cases, the general trend is the same as that in the sloshing test case. The movement according to the change in streamlines is the most accurate, followed by the second order method; while the first order method is the least accurate. For the larger time step considered, the errors in the streamline and first order methods are significantly larger than the other two methods, while the difference is not as big for the smaller time step. The error in the first order method with the large time step is too huge for this combination to be used in any practical simulations. The relative errors in volume at $t=3 s$ for each of the methods, and for both time steps used is also tabulated in Table~\ref{tab:TankFillingTable}.
\begin{table}
\caption{Tank filling test case: Relative errors in numerical volume when compared with the analytical value, for different methods of movement at $t=3 s$. The evolution of these errors is shown in Figure~\ref{Fig:TankFillingResults}. M1 stands for the first order method of movement, M2 the second order method, M3 the movement along the streamline, and M4 the movement according to the change of streamlines.}
\centering
\label{tab:TankFillingTable}
	\begin{tabular}{l|c|c}
\hline
&& \\[\dimexpr-\normalbaselineskip+2pt]
	 & $\Delta t = 0.6\times 10^{-2}$ & $\Delta t = 0.6\times 10^{-3}$\\
\hline
&& \\[\dimexpr-\normalbaselineskip+2pt]
	M1 & $60.4\%$ & $8.4\%$ \\
	M2 & $10.8\%$ & $4.1\%$ \\
	M3 & $29.5\%$ & $6.0\%$ \\
	M4 & $5.5\%$ & $2.6\%$ \\			
\hline	
	\end{tabular}
\end{table}

\section{Conclusion}
\label{sec:Conclusion}

We presented different methods of moving point clouds for fully Lagrangian meshfree methods. The most often used first-order method was shown to produce very inaccurate results. We presented two different methods to improve accuracy of movement. In the first method, each point is moved along approximated streamlines in the neighbourhood of that point. In the second method, the change of streamlines around that point from the previous time step to the current time step is considered, and movement is done based on that. Each of these methods provide a better way to track particle motion, and thus they improve the accuracy of the Lagrangian movement step.

The numerical results match the theoretical expectations. For quasi-stationary flow, moving points along the velocity streamlines produces better results than the other methods. On the other hand, for situations with rapidly changing flow profiles, movement according to the change of streamlines showed the best results.

A key point of this work is to emphasize that the accuracy and conservative properties of fully Lagrangian meshfree methods can be improved without the need of resorting to artificial, non-physical corrections of particle locations or the physical quantities carried by the particles.

%

\end{document}

%% file: Pat.pdf_tex

\begingroup
  \makeatletter
  \providecommand\color[2][]{%
    \errmessage{(Inkscape) Color is used for the text in Inkscape, but the package 'color.sty' is not loaded}
    \renewcommand\color[2][]{}%
  }
  \providecommand\transparent[1]{%
    \errmessage{(Inkscape) Transparency is used (non-zero) for the text in Inkscape, but the package 'transparent.sty' is not loaded}
    \renewcommand\transparent[1]{}%
  }
  \providecommand\rotatebox[2]{#2}
  \ifx\svgwidth\undefined
    \setlength{\unitlength}{479.93245869pt}
  \else
    \setlength{\unitlength}{\svgwidth}
  \fi
  \global\let\svgwidth\undefined
  \makeatother
  \begin{picture}(1,0.7468464)%
    \put(0,0){\includegraphics[width=\unitlength]{./Pat.pdf}}%
    \put(0.32957111,0.1613148){\color[rgb]{0,0,0}\makebox(0,0)[rb]{\smash{$t^{n-1}$}}}%
    \put(0.53793375,0.52803306){\color[rgb]{0,0,0}\makebox(0,0)[rb]{\smash{$t^{n}$}}}%
    \put(0.90465202,0.69472318){\color[rgb]{0,0,0}\makebox(0,0)[rb]{\smash{$\vec{x}(t^{n+1}) = ?$}}}%
    \put(0.41291617,0.32800492){\color[rgb]{0,0,0}\makebox(0,0)[rb]{\smash{$t^{n-1}+\tau$}}}%
    \put(0.70462387,0.61137812){\color[rgb]{0,0,0}\makebox(0,0)[rb]{\smash{$t^{n}+\tau$}}}%
  \end{picture}%
\endgroup